\documentclass{article}
\usepackage{amsmath}
\usepackage{amssymb}

\newcounter{theorem} 
\newcounter{lemma} 

\renewcommand{\thetheorem}{\arabic{theorem}}
\renewcommand{\thelemma}{\arabic{lemma}}
\newcommand{\theor}{\par\refstepcounter{theorem}%
 {\bf Theorem \thetheorem .}\,\,}

\newcommand{\lem}{\par\refstepcounter{lemma}%
{\bf Lemma \thelemma .}\,\,}
\,
\,

\begin{document}
\Large
\begin{center}
\textbf{Curvilinear integral theorems for monogenic\vskip2mm
 functions in
commutative associative algebras}\end{center}\large
\vskip6mm\begin{center}
 \textbf{V.~S.~Shpakivskyi}\end{center}

\vskip5mm

\textbf{Abstract.} We consider an arbitrary finite-dimensional
commutative associative algebra, $\mathbb{A}_n^m$, with unit over the field of complex number
 with $m$ idempotents. Let
$e_1=1,e_2,e_3$ be elements of $\mathbb{A}_n^m$ which are linearly
independent over the field of real numbers. We consider monogenic
(i.~e. continuous and differentiable in the sense of Gateaux)
functions of the variable $xe_1+ye_2+ze_3$\,, where $x,y,z$ are
real. For mentioned monogenic function we prove curvilinear analogues of the
Cauchy integral theorem, the Morera theorem and the Cauchy integral formula.

\vskip5mm
\textbf{Keywords:} Commutative associative algebra; Cauchy integral theorem; Morera theorem;
Cauchy integral formula.
\vskip5mm

\section{Introduction.}
\vskip2mm

The Cauchy integral theorem and Cauchy integral
formula for the holomorphic function of the complex variable are
a fundamental result of the classical complex analysis.
Analogues of these results are also an important tool in commutative
algebras of dimensional more that $2$.

In the paper of E.~R.~Lorch \cite{Lorch} for functions differentiable in the
sense of Lorch in an arbitrary convex domain of commutative
associative Banach algebra, some properties similar to properties
of holomorphic functions of complex variable (in particular, the
curvilinear integral Cauchy theorem and the integral Cauchy formula, the
Taylor expansion and the Morera theorem) are established. E.~K.~Blum \cite{Blum}
withdrew a convexity condition of a domain in the
mentioned results from \cite{Lorch}.

Let us note that {\it a
priori}\/ the differentiability of a function in the
sense of Gateaux is a restriction weaker than the
differentiability of this function in the sense of Lorch.
Therefore, we consider a \textit{monogenic} functions defined as a continuous
and differentiable in the sense of Gateaux. Also we assume that a
monogenic function is given in a domain of three-dimensional subspace of
an arbitrary commutative associative algebra with unit over the field of complex
 numbers.
 In this situation
the results established in the papers \cite{Lorch,Blum} is not applicable
 for a mentioned monogenic function, because it deals with an
 integration along a curve on which the function is not given, generally speaking.

In the papers \cite{Pl-Shp3,Pl-Shp-Algeria,Pl-Pukh-Analele} for monogenic
function the curvilinear analogues of the
Cauchy integral theorem, the Cauchy integral formula and the
 Morera theorem are obtained in special finite-dimensional
 commutative associative algebras.

 In this paper we generalize results of the papers
 \cite{Pl-Shp3,Pl-Shp-Algeria,Pl-Pukh-Analele} for
an arbitrary commutative associative algebra over the field of
complex numbers.

 Let us note that some analogues of the curvilinear  Cauchy integral theorem
 and the Cauchy integral formula for another classes of functions in
 special commutative algebras are
established in the papers \cite{Goncharow,Ketchum-28,Ketchum-29,Rosculet-54,Rosculet-55}.

\section{The algebra $\mathbb{A}_n^m$.}
\vskip2mm

Let $\mathbb{N}$ be the set of natural numbers.
We fix the numbers $m,n\in\mathbb{N}$ such that $m\leq n$.
Let $\mathbb{A}_n^m$ be an arbitrary commutative associative algebra with
 unit over the field of complex number $\mathbb{C}$.
   E.~Cartan \cite[pp.~33 -- 34]{Cartan}
  proved that in the algebra $\mathbb{A}_n^m$ there exist a basis $\{I_k\}_{k=1}^{n}$
  satisfies the following multiplication rules:
  \vskip3mm
1.  \,\,  $\forall$\; $r,s\in[1,m]\cap\mathbb{N}$\,: \qquad $I_rI_s=\left\{
\begin{array}{rcl}
0 &\mbox{if} & r\neq s,\vspace*{2mm} \\
I_r &\mbox{if} & r=s;\\
\end{array}
\right.$

\vskip5mm

2. \,\,  $\forall$\; $r,s\in[m+1,n]\cap\mathbb{N}$\,: \qquad $I_rI_s=
\sum\limits_{k=\max\{r,s\}+1}^n\Upsilon_{r,k}^{s}I_k$\,;

\vskip5mm

3.\,\, $\forall$\; $s\in[m+1,n]\cap\mathbb{N}$\;  $\exists!\;
 u_s\in[1,m]\cap\mathbb{N}$ \;$\forall$\,
 $r\in[1,m]\cap\mathbb{N}$\,:\;\;

\begin{equation}\label{mult_rule_3}
I_rI_s=\left\{
\begin{array}{ccl}
0 \;\;\mbox{if}\;\;  r\neq u_s\,,\vspace*{2mm}\\
I_s\;\;\mbox{if}\;\;  r= u_s\,. \\
\end{array}
\right.\medskip
\end{equation}
Furthermore, the structure constants $\Upsilon_{r,k}^{s}\in\mathbb{C}$
 satisfy the associativity conditions:
\vskip2mm
(A\,1).\,\, $(I_rI_s)I_p=I_r(I_sI_p)$ \; $\forall$\, $r,s,p\in[m+1,n]\cap\mathbb{N}$;
\vskip2mm
(A\,2).\,\, $(I_uI_s)I_p=I_u(I_sI_p)$ \; $\forall$\, $u\in[1,m]\cap\mathbb{N}$\;
 $\forall$\,
 $s,p\in[m+1,n]\cap\mathbb{N}$.
\vskip2mm
Obviously, the first $m$ basis vectors $\{I_u\}_{u=1}^m$ are the
idempotents and, respectively,
form the semi-simple subalgebra. Also the vectors $\{I_r\}_{r=m+1}^n$
form the nilpotent subalgebra of algebra
 $\mathbb{A}_n^m$.
The unit of $\mathbb{A}_n^m$ is the element $1=\sum_{u=1}^mI_u$. Therefore,
we will write that the algebra $\mathbb{A}_n^m$ is a semi-direct sum of the
$m$-dimensional semi-simple subalgebra $S$ and $(n-m)$-dimensional
nilpotent subalgebra $N$, i.~e.
$$\mathbb{A}_n^m=S\oplus_s N.
$$

In the cases where $\mathbb{A}_n^m$ has some specific properties,
the following propositions are true.

\textbf{Proposition 1 \cite{Shpakivskyi-2014}.} \textit{If there exists the unique
$u_0\in[1,m]\cap\mathbb{N}$ such that $I_{u_0}I_s=I_s$ for all
$s=m+1,\ldots,n$,  then the associativity condition \em (A\,2) \em
is satisfied.}

Thus, under the conditions of Proposition 1,
the associativity condition (A\,1) is only required.
It means that the nilpotent subalgebra of $\mathbb{A}_n^m$ with
the basis $\{I_r\}_{r=m+1}^n$  can be an arbitrary commutative
associative nilpotent algebra of dimension $n-m$. We note that such
nilpotent algebras are fully described for the dimensions
$1,2,3$ in the paper \cite{Burde_de_Graaf}, and some four-dimensional nilpotent algebras
can be found in the papers \cite{Burde_Fialowski},  \cite{Martin}.

\textbf{Proposition 2 \cite{Shpakivskyi-2014}.} \textit{If all $u_r$ are different in the
multiplication rule \em 3\em ,
 then $I_sI_p=0$ for all $s,p=m+1,\ldots, n$.}

Thus, under the conditions of Proposition 2,
the multiplication table of the
nilpotent subalgebra of $\mathbb{A}_n^m$ with the basis
$\{I_r\}_{r=m+1}^n$ consists only of zeros, and all associativity
conditions are satisfied.

The algebra $\mathbb{A}_n^m$ contains $m$ maximal ideals
$$\mathcal{I}_u:=\Biggr\{\sum\limits_{k=1,\,k\neq u}^n\lambda_kI_k:\lambda_k\in
\mathbb{C}\Biggr\}, \quad  u=1,2,\ldots,m,
$$
the intersection of which is the radical $$\mathcal{R}:=
\Bigr\{\sum\limits_{k=m+1}^n\lambda_kI_k:\lambda_k\in
\mathbb{C}\Bigr\}.$$

We define $m$ linear functionals $f_u:\mathbb{A}_n^m\rightarrow\mathbb{C}$ by put
$$f_u(I_u)=1,\quad f_u(\omega)=0\quad\forall\,\omega\in\mathcal{I}_u\,,
\quad u=1,2,\ldots,m.
$$
Since the kernels of functionals $f_u$ are, respectively, the maximal ideals
 $\mathcal{I}_u$, then these functionals are also continuous and multiplicative
  (see \cite[p. 147]{Hil_Filips}).

\vskip2mm
\section{Monogenic functions.}
\vskip2mm

We consider the vectors $e_1=1,e_2,e_3$ in $\mathbb{A}_n^m$ which are linearly
independent over the field of real number $\mathbb{R}$. It means that the equality
$$\alpha_1e_1+\alpha_2e_2+\alpha_3e_3=0,\quad \alpha_1,\alpha_2,
\alpha_3\in\mathbb{R},$$
holds if and only if $\alpha_1=\alpha_2=
\alpha_3=0$.

Let the vectors $e_1=1,e_2,e_3$ have the following decompositions with respect to
the basis $\{I_k\}_{k=1}^n$:
\begin{equation}\label{e_1_e_2_e_3}
e_1=1,\quad e_2=\sum\limits_{k=1}^na_kI_k\,,\quad e_3=\sum\limits_{k=1}^nb_kI_k\,,
\end{equation}
where $a_k,b_k\in\mathbb{C}$.

Let $\zeta:=xe_1+ye_2+ze_3$, where $x,y,z\in\mathbb{R}$.
It is also obvious that
 $\xi_u:=f_u(\zeta)=x+ya_u+zb_u$,\, $u=1,2,\ldots,m$.
 Let
 $E_3:=\{\zeta=xe_1+ye_2+ze_3:\,\, x,y,z\in\mathbb{R}\}$ be the
linear span of vectors $e_1,e_2,e_3$ over the field  of real numbers
$\mathbb{R}$. We note that in the further investigations,
 it is essential assumption:  $f_u(E_3)=\mathbb{C}$ for all $u=1,2,\ldots,m$,
 where $f_u(E_3)$ is the image of $E_3$ under the mapping $f_u$.
 Obviously, it holds if and only if for every fixed $u=1,2, \ldots, m$
at least one of the numbers $a_u$ or $b_u$ belongs to $\mathbb{C}\setminus\mathbb{R}$.

With a set $Q\subset\mathbb{R}^3$ we associate the set $Q_\zeta:=
\{\zeta=xe_1+ye_2+ze_3:(x,y,z)\in Q\}$ in $E_3$. We also note that
the topological properties of a set $Q_\zeta$ in $E_3$ understood as a
corresponding topological properties of a set $Q$ in $\mathbb{R}^3$.
For example, a homotopicity of a curve $\gamma_\zeta\subset E_3$ to the zero means a
homotopicity of $\gamma\subset\mathbb{R}^3$ to the zero, etc.

Let $\Omega$ be a domain in ${\mathbb R}^3$.

A continuous function
$\Phi:\Omega_{\zeta}\rightarrow\mathbb{A}_n^m$ is \textit{monogenic}
in $\Omega_{\zeta}$ if $\Phi$ is differentiable in the sense of
Gateaux in every point of $\Omega_{\zeta}$, i.~e. if  for every
$\zeta\in\Omega_{\zeta}$ there exists an element
$\Phi'(\zeta)\in\mathbb{A}_n^m$ such that
\begin{equation}\label{monogennaOZNA}\medskip
\lim\limits_{\varepsilon\rightarrow 0+0}
\left(\Phi(\zeta+\varepsilon
h)-\Phi(\zeta)\right)\varepsilon^{-1}= h\Phi'(\zeta)\quad\forall\,
h\in E_{3}.\medskip
\end{equation}
$\Phi'(\zeta)$ is the \textit{Gateaux derivative} of the function
$\Phi$ in the point $\zeta$.

Consider the decomposition of a function
$\Phi:\Omega_{\zeta}\rightarrow\mathbb{A}_n^m$ with respect to the
basis $\{I_k\}_{k=1}^n$:
\begin{equation}\label{rozklad-Phi-v-bazysi}
\Phi(\zeta)=\sum_{k=1}^n U_k(x,y,z)\,I_k\,.
 \end{equation}

In the case where the functions $U_k:\Omega\rightarrow\mathbb{C}$ are
$\mathbb{R}$-differentiable in $\Omega$, i.~e. for every
$(x,y,z)\in\Omega$,
$$U_k(x+\Delta x,y+\Delta y,z+\Delta z)-U_k(x,y,z)=
\frac{\partial U_k}{\partial x}\,\Delta x+
\frac{\partial U_k}{\partial y}\,\Delta y+\frac{\partial
U_k}{\partial z}\,\Delta z+$$ $$+\,o\left(\sqrt{(\Delta
x)^2+(\Delta y)^2+(\Delta z)^2}\,\right), \qquad (\Delta
x)^2+(\Delta y)^2+(\Delta z)^2\to 0\,,$$
the function $\Phi$ is monogenic in the domain $\Omega_{\zeta}$ if
and only if the following Cauchy~-- Riemann conditions are
satisfied in $\Omega_{\zeta}$:
\begin{equation}\label{Umovy_K-R}
\frac{\partial \Phi}{\partial y}=\frac{\partial \Phi}{\partial
x}\,e_{2},\quad \frac{\partial \Phi}{\partial z}=\frac{\partial
\Phi}{\partial x}\,e_{3}.
\end{equation}

Expansion of the resolvent is of the form
\begin{equation}\label{rozkl-rezol-A_n^m}
(te_1-\zeta)^{-1}=\sum\limits_{u=1}^m\frac{1}{t-\xi_u}\,I_u+
 \sum\limits_{s=m+1}^{n}\sum\limits_{k=2}^{s-m+1}\frac{Q_{k,s}}
 {\left(t-\xi_{u_{s}}\right)^k}\,I_{s}\,
  \end{equation}
  $$ \forall\,t\in\mathbb{C}:\,
t\neq \xi_u,\quad u=1,2,\ldots,m,$$
where $Q_{k,s}$ are determined by the following recurrence
relations:
\begin{equation}\label{lem_3_2-}
Q_{2,s}:=T_{s}\,,\quad
Q_{k,s}=\sum\limits_{r=k+m-2}^{s-1}Q_{k-1,r}\,B_{r,\,s}\,,\; \;\;k=3,4,\ldots,s-m+1.
\end{equation}
with
$$T_s:=ya_s+zb_s\,,\;\;B_{r,s}:=\sum\limits_{k=m+1}^{s-1}T_k
\Upsilon_{r,s}^k\,,\; \;\;s=m+2,\ldots,n,$$
 and natural numbers $u_s$ are defined in the
rule  3  of the multiplication table of the algebra $\mathbb{A}_n^m$.

From the relations  (\ref{rozkl-rezol-A_n^m}) follows that the points
 $(x,y,z)\in\mathbb{R}^3$
corresponding to the noninvertible elements
 $\zeta\in\mathbb{A}_n^m$ form the straight lines
  \[L_u:\quad\left\{
\begin{array}{r}x+y\,{\rm Re}\,a_u+z\,{\rm Re}\,b_u=0,\vspace*{3mm} \\
y\,{\rm Im}\,a_u+z\,{\rm Im}\,b_u=0 \\ \medskip
\end{array} \right.\]
in the three-dimensional space $\mathbb{R}^3$.

Denote by $D_u\subset\mathbb{C}$ the image of $\Omega_\zeta$ under the mapping
$f_u$,\, $u=1,2, \ldots, m$.
A constructive description of all monogenic functions in the algebra
$\mathbb{A}_n^m$
 by means of holomorphic functions of the complex variable are
  obtained in the paper \cite{Shpakivskyi-2014}.  Namely, it is
 proved the theorem:

Let a domain $\Omega\subset
\mathbb{R}^{3}$ be convex in the
direction of the straight lines $L_u$ and $f_u(E_3)=\mathbb{C}$ for all
 $u=1,2,\ldots, m$.  Then any
monogenic function $\Phi:\Omega_{\zeta}\rightarrow\mathbb{A}_n^m$
can be expressed in the form
 \begin{equation}\label{Teor--1}
\Phi(\zeta)=\sum\limits_{u=1}^mI_u\,\frac{1}{2\pi i}\int\limits_{\Gamma_u}
F_u(t)(te_1-\zeta)^{-1}\,dt+
\sum\limits_{s=m+1}^nI_s\,\frac{1}{2\pi i}\int\limits_
{\Gamma_{u_s}}G_s(t)(te_1-\zeta)^{-1}\,dt,
 \end{equation}
where $F_u$ is the certain holomorphic function in
a domain  $D_u$;  $G_s$  is the certain holomorphic function in
a domain
 $D_{u_s}$; $\Gamma_q$ is a closed Jordan rectifiable curve
lying in the domain $D_q$ surround
 a point $\xi_q$ and containing no points
$\xi_{\ell}$, $\ell,q=1,2,\ldots, m$,\,$\ell\neq q$.

\vskip3mm
\section{Cauchy integral theorem for a curvilinear integral.}
\vskip3mm

 Let $\gamma$ be a Jordan rectifiable curve in $\mathbb{R}^{3}$. For a
continuous function $\Psi:\gamma_{\zeta}\rightarrow
\mathbb{A}_n^m$ of the form
\begin{equation}\label{Phi-form}
\Psi(\zeta)=\sum\limits_{k=1}^{n}{U_k(x,y,z)\,I_k}+i\sum
\limits_{k=1}^{n}{V_k(x,y,z)\,I_k},
\end{equation}
where $(x,y,z)\in\gamma$ and $U_k:\gamma\rightarrow\mathbb{R}$,
$V_k:\gamma\rightarrow\mathbb{R}$,
we define an integral along a Jordan rectifiable curve $\gamma_\zeta$ by
 the equality:
$$\int\limits_{\gamma_{\zeta}}\Psi(\zeta)d\zeta:=\sum\limits_{k=1}^{n}
I_{k}\int\limits_{\gamma}U_{k}(x,y,z)dx+
\sum\limits_{k=1}^ne_{2}I_{k}\int\limits_{\gamma}U_{k}(x,y,z)dy+$$
$$+\sum\limits_{k=1}^ne_{3}I_{k}\int\limits_{\gamma}U_{k}(x,y,z)dz
+i\sum\limits_{k=1}^nI_{k}\int\limits_{\gamma}V_{k}(x,y,z)dx+$$
$$+i\sum\limits_{k=1}^ne_{2}I_{k}\int\limits_{\gamma}V_{k}(x,y,z)dy+
i\sum\limits_{k=1}^ne_{3}I_{k}\int\limits_{\gamma}V_{k}(x,y,z)dz,$$
where $d\zeta:=dx+e_{2}dy+e_{3}dz$.

Also we define a surface integral. Let $\Sigma$ be a piece-smooth surface in
 $\mathbb{R}^{3}$. For a continuous
function $\Psi:\Sigma_{\zeta}\rightarrow \mathbb{A}_n^m$ of the
form (\ref{Phi-form}), where $(x,y,z)\in\Sigma$ and $U_k:\Sigma\rightarrow\mathbb{R}$,
$V_k:\Sigma\rightarrow\mathbb{R}$, we define a surface
integral on $\Sigma_{\zeta}$ with the differential form
$dxdy$, by the equality
$$\int\limits_{\Sigma_{\zeta}}\Psi(\zeta)dxdy:=
\sum\limits_{k=1}^nI_{k}\int\limits_{\Sigma}U_{k}(x,y,z)dxdy+
i\sum\limits_{k=1}^nI_{k}\int\limits_{\Sigma}V_{k}(x,y,z)dxdy.
$$
A similarly defined the integrals with the forms $dydz$ and $dzdx$.

If a function $\Phi:\Omega_\zeta\rightarrow\mathbb{A}_n^m$ is
continuous together with partial derivatives of the first order in
a domain $\Omega_\zeta$, and $\Sigma$ is a piece-smooth surface
in $\Omega$, and the edge $\gamma$ of surface $\Sigma$ is a
rectifiable Jordan curve, then the following analogue of the
Stokes formula is true:
$$\int\limits_{\gamma_{\zeta}}\Psi(\zeta)d\zeta=\int\limits_
{\Sigma_{\zeta}}\left(\frac{\partial\Psi}{\partial
x}e_{2}-\frac{\partial\Psi}{\partial y}\right)dxdy
+\left(\frac{\partial\Psi}{\partial
y}e_{3}-\frac{\partial\Psi}{\partial z}e_{2}\right)dydz+$$
\begin{equation}\label{form-Stoksa}+\left(\frac{\partial\Psi}{\partial
z}-\frac{\partial\Psi}{\partial x}e_{3}\right)dzdx.
\end{equation}
Now, the next theorem is a result of the formula
(\ref{form-Stoksa}) and the equalities (\ref{Umovy_K-R}).
\vskip2mm

\theor\label{teo-int-po-kryv-z-neper-poh}
\emph{ Suppose that $\Phi:\Omega_{\zeta}\rightarrow\mathbb{A}_n^m$ is a
monogenic function in a domain $\Omega_{\zeta}$, and $\Sigma$ is a
piece-smooth surface in $\Omega$, and the edge $\gamma$ of surface
$\Sigma$ is a rectifiable Jordan curve. Then}
\begin{equation}\label{form-Koshi-po-kryv}
\int\limits_{\gamma_{\zeta}}\Phi(\zeta)d\zeta=0.
\end{equation}

In the case where a domain $\Omega$ is convex, then by the usual way
 (see, e.~g., \cite{Privalov}) the equality (\ref{form-Koshi-po-kryv})
  can be prove for an arbitrary closed Jordan rectifiable curve $\gamma_\zeta$.

In the case where a domain $\Omega$ is an arbitrary, then similarly
 to the proof of Theorem 3.2 \cite{Blum} we can
prove the following

\theor\label{teo-int-po-kryv-Blum}
 \emph{ Let $\Phi:\Omega_{\zeta}\rightarrow\mathbb{A}_n^m$ be a monogenic
function in a domain $\Omega_{\zeta}$. Then for every closed
Jordan rectifiable curve $\gamma$ homotopic to a point in
$\Omega$, the equality \em (\ref{form-Koshi-po-kryv}) \em is true.}

\vskip3mm
\section{The Morera theorem.}
\vskip3mm

To prove the analogue of Morera theorem in the algebra $\mathbb{A}_n^m$,
we introduce auxiliary notions and prove some auxiliary statements.

Let us consider the algebra $\mathbb{A}_n^m(\mathbb{R})$ with  the
 basis $\{I_k,iI_k\}_{k=1}^n$ over the field $\mathbb{R}$ which is isomorphic to the
algebra $\mathbb{A}_n^m$ over the field $\mathbb{C}$. In the algebra
$\mathbb{A}_n^m(\mathbb{R})$ there exist another basis $\{e_k\}_{k=1}^{2n}$,
where the vectors $e_1,e_2,e_3$ are the same as in the Section 3.

For the element $a:=\sum\limits_{k=1}^{2n}a_ke_k$,\, $a_k\in\mathbb{R}$ we define
the Euclidian norm $$\|a\|:=\sqrt{\sum\limits_{k=1}^{2n}a_k^2}\,.$$ Accordingly,
$\|\zeta\|=\sqrt{x^2+y^2+z^2}$ and $\|e_1\|=\|e_2\|=\|e_3\|=1$.

Using the Theorem on equivalents of norms, for the element
$b:=\sum\limits_{k=1}^{n}(b_{1k}+ib_{2k})I_k$,\, $b_{1k},b_{2k}\in\mathbb{R}$ we have
the following inequalities
\begin{equation}\label{ner-dlja-integrala-dop}
|b_{1k}+ib_{2k}|\leq\sqrt{\sum\limits_{k=1}^{2n}\big(b_{1k}^2+b_{2k}^2\big
)}\,\leq c \|b\|,
\end{equation}
where $c$ is a positive constant does not depend on $b$.
\vskip2mm

\lem\label{lem-22} \textit{If $\gamma$ is a closed Jordan
rectifiable curve in $\mathbb{R}^{3}$ and function
$\Psi:\gamma_{\zeta}\rightarrow
\mathbb{A}_n^m$ is continuous, then}

\begin{equation}\label{ner-dlja-integrala}
\Biggr\|\int\limits_{\gamma_{\zeta}}\Psi(\zeta)\,d\zeta\Biggr\|\leq
c \int\limits_{\gamma_{\zeta}}\|\Psi(\zeta)\|\|d\zeta\|,
\end{equation}
\textit{where
$c$ is a positive absolutely constant.}

\vskip1mm
\textbf{Proof.} Using the representation of function $\Psi$ in the form
 (\ref{Phi-form}) for $(x,y,z)\in\gamma$, we obtain \vspace{2mm}
 $$\Biggl\|\int\limits_{\gamma_{\zeta}}\Psi(\zeta)d\zeta\Biggr\|\le
 \sum\limits_{k=1}^{n}\|I_k\|\int\limits_{\gamma}\bigl|U_{k}(x,y,z)+iV_{k}(x,y,z)\bigr|
 \,dx+$$\vspace{1mm}
$$+\sum\limits_{k=1}^{n}\|e_{2}I_{k}\|\int\limits_{\gamma}\bigl|U_{k}(x,y,z)+iV_{k}
(x,y,z)\bigr|\,dy+$$\vspace{1mm}
$$+\sum\limits_{k=1}^{n}\|e_{3}I_{k}\|\int\limits_{\gamma}\bigl|U_{k}(x,y,z)+iV_{k}
(x,y,z)\bigr|\,dz.$$
Now, taking into account the inequality (\ref{ner-dlja-integrala-dop}) for
 $b=\Psi(\zeta)$ and the inequalities $\|e_sI_k\|\leq c_s$,\, $s=1,2,3$,
 where
$c_s$ are positive absolutely constants,
 we obtain the relation (\ref{ner-dlja-integrala}).
\noindent The lemma is proved.

Using Lemma \ref{lem-22}, for functions taking values in the algebra $\mathbb{A}_n^m$, the
following Morera theorem can be established in the usual way.
\vskip2mm

\theor\label{teo_Morera}
 \emph{ If a function $\Phi:\Omega_{\zeta}\rightarrow\mathbb{A}_n^m$ is
continuous in a domain $\Omega_{\zeta}$ and satisfies the equality
\begin{equation} \label{Morera}
\int\limits_{\partial\triangle_\zeta}\Phi(\zeta)d\zeta=0
\end{equation}
 for every triangle $\triangle_\zeta$ such that closure
$\overline{\triangle_\zeta}\subset\Omega_\zeta$, then the function
$\Phi$ is monogenic in the domain $\Omega_{\zeta}$.}

\vskip3mm
\section{Cauchy integral formula for a curvilinear integral.}
\vskip3mm

 Let $\zeta_0:=x_0e_1+y_0e_2+z_0e_3$ be a point in a domain
$\Omega_\zeta\subset E_3$. In a neighborhood of $\zeta_0$ contained in
$\Omega_\zeta$ let us take a circle $C_\zeta(\zeta_0,\varepsilon)$ of radius
$\varepsilon$ with
the center at the point
$\zeta_0$. By $C_u(\xi^{(0)}_u,\varepsilon)\subset\mathbb{C}$ we denote the image of $C_\zeta(\zeta_0,\varepsilon)$
under the mapping $f_u$, $u=1,2,\ldots,m$.  We assume that
the circle $C_\zeta(\zeta_0,\varepsilon)$ \emph{embraces the set} $\{\zeta-\zeta_0:
(x,y,z)\in\bigcup\limits_{u=1}^m L_u\}$.
It means that the curve $C_u(\xi^{(0)}_u,\varepsilon)$ bounds some domain $D_u'$ and
$f_u(\zeta_0)=\xi^{(0)}_{u}\in D_u'$, \, $u=1,2,\ldots,m$.

We say that the curve $\gamma_\zeta\subset\Omega_\zeta$ \emph{embraces once the set} $\{\zeta-\zeta_0:(x,y,z)\in\bigcup\limits_{u=1}^m L_u\}$, if there exists a
 circle $C_\zeta(\zeta_0,\varepsilon)$
which embraces the mentioned set and is homotopic to $\gamma_\zeta$ in the domain $\Omega_\zeta\setminus\{\zeta-\zeta_0:(x,y,z)\in\bigcup\limits_{u=1}^m L_u\}$.

Since the function $\zeta^{-1}$ is continuous on the curve $C_\zeta(0,\varepsilon)$, then there exist the integral
\begin{equation}\label{lambda}
\lambda:=\int\limits_{C_\zeta(0,\varepsilon)}\zeta^{-1}d\zeta.
\end{equation}

The following theorem is an analogue of Cauchy integral
theorem for monogenic function
 $\Phi:\Omega_{\zeta}\rightarrow\mathbb{A}_{n}^m$.\vskip2mm

\theor\label{teo-formula-Koshi} \textit{Suppose that a domain $\Omega\subset
\mathbb{R}^{3}$ is convex in the
direction of the straight lines $L_u$ and $f_u(E_3)=\mathbb{C}$ for all
 $u=1,2,\ldots, m$. Suppose also that
 $\Phi:\Omega_{\zeta}\rightarrow\mathbb{A}_n^m$ is a monogenic
function in $\Omega_\zeta$.
 Then for every point
$\zeta_{0}\in\Omega_{\zeta}$ the following equality is true:
\begin{equation}\label{form-Koshi}
\lambda\,\Phi(\zeta_{0})=
\int\limits_{\gamma_{\zeta}}\Phi(\zeta)\left(\zeta-\zeta_{0}\right)^{-1}d\zeta,
\end{equation}
where $\gamma_{\zeta}$ is an arbitrary closed Jordan rectifiable curve in
$\Omega_\zeta$, that embraces once the set $\{\zeta-\zeta_0:(x,y,z)\in
\bigcup\limits_{u=1}^m L_u\}$.
}\vskip 1mm

\textbf{Proof.} Inasmuch as $\gamma_\zeta$ is homotopic to $C_\zeta(\zeta_0,\varepsilon)$ in the
domain $\Omega_\zeta\setminus\{\zeta-\zeta_0:(x,y,z)\in\bigcup\limits_{u=1}^m
 L_u\}$, it follows from Theorem \ref{teo-int-po-kryv-Blum} that
   \begin{equation}\label{form-Koshi-}
\int\limits_{\gamma_\zeta}\Phi(\zeta)\left(\zeta-\zeta_{0}\right)^{-1}d\zeta=
\int\limits_{C_\zeta(\zeta_0,\varepsilon)}\Phi(\zeta)\left(\zeta-\zeta_{0}\right)
^{-1}d\zeta.
\end{equation}

Taking into account the equality (\ref{form-Koshi-}) we represent the
integral on the right-hand side of
equality (\ref{form-Koshi}) as the sum of the following two
integrals:\vspace{-0.5mm}
$$\int\limits_{\gamma_{\zeta}}\Phi(\zeta)\left(\zeta-\zeta_{0}\right)^{-1}d\zeta=
\int\limits_{C_\zeta(\zeta_0,\varepsilon)}
(\Phi(\zeta)-\Phi(\zeta_{0}))\left(\zeta-\zeta_{0}\right)^{-1}d\zeta+$$\vspace{-3mm}
$$+\Phi(\zeta_{0})\int\limits_{C_\zeta(\zeta_0,\varepsilon)}\left(\zeta-\zeta_{0}\right)^{-1}d\zeta=
:J_{1}+J_{2}.$$

Let us note that from the relation (\ref{form-Koshi-}) follows that if there exist
the integral in the equality (\ref{lambda}) then it does not depend on
$\varepsilon$. As a consequence of the equalities (\ref{lambda}), (\ref{form-Koshi-}),
we have the following relation
\begin{equation}\label{rav-2-}
J_2=
\Phi(\zeta_{0})\int\limits_{C_\zeta(0,\varepsilon)}\tau^{-1}d\tau=\lambda\,
\Phi(\zeta_{0}),
\end{equation}
where $\tau:=\zeta-\zeta_0$.

The integrand in the integral $J_1$ is bounded by a constant which
does not depend on $\varepsilon$: when $\varepsilon\rightarrow0$
 the integrand tends to $\Phi'(\zeta_0)$ (see Lemma 5 \cite{Shpakivskyi-2014}).
Therefore, using the Lemma \ref{lem-22} the integral $J_1$ tends to zero
as $\varepsilon\rightarrow0$.
The theorem is proved.

Below, it will be shown that the constant $\lambda$ is an invertible element in $\mathbb{A}_n^m$.

\vskip3mm
\section{A constant $\lambda$.}
\vskip3mm

In some special algebras (see \cite{Pl-Shp3,Pl-Shp-Algeria,Pl-Pukh-Analele})
the Cauchy integral formula (\ref{form-Koshi}) has the form
\begin{equation}\label{form-Koshi--}
\Phi(\zeta_{0})=\frac{1}{2\pi i}
\int\limits_{\gamma_{\zeta}}\Phi(\zeta)\left(\zeta-\zeta_{0}\right)^{-1}d\zeta,
\end{equation}
i.~e.
\begin{equation}\label{lambda-}
\lambda=2\pi i.
\end{equation}

In this Section we indicate a set of algebras $\mathbb{A}_n^m$ for which
(\ref{lambda-}) holds. In this a way we first consider some auxiliary statements.

As a consequence of the expansion (\ref{rozkl-rezol-A_n^m}), we obtain the
following equality:
\begin{equation}\label{obr-elem}
\zeta^{-1}=\sum\limits_{k=1}^n\widetilde{A}_k\,I_k
\end{equation}
with the coefficients\, $\widetilde{A}_k$\, determined by the following
 relations:

\begin{equation}\label{A__p}
\begin{array}{c}
\displaystyle
\widetilde{A}_u=\frac{1}{\xi_u}\,,\;\;u=1,2,\ldots,m, \quad
 \vspace*{4mm}\\
 \displaystyle
\widetilde{A}_s=\sum\limits_{k=2}^{s-m+1}\frac{\widetilde{Q}_{k,s}}{\xi_{u_s}^k}\,,\quad s=m+1,m+2,\ldots,n,
\end{array}
\end{equation}
where $\widetilde{Q}_{k,s}$ are determined by the following recurrence
relations:
\begin{equation}\label{lem_3_2--}
\widetilde{Q}_{2,s}:=-T_{s}\,,\quad
\widetilde{Q}_{k,s}=-\sum\limits_{r=k+m-2}^{s-1}\widetilde{Q}_{k-1,r}\,B_{r,\,s}\,,\; \;\;k=3,4,\ldots,s-m+1.
\end{equation}
where $T_s$ and $B_{r,s}$ are the same as in the equalities (\ref{lem_3_2-}),
 and natural numbers $u_s$ are defined in the
rule  3  of the multiplication table of the algebra $\mathbb{A}_n^m$.

Taking into account the equality
(\ref{obr-elem}) and the relation
$$d\zeta=dxe_1+dye_2+dze_3=\sum\limits_{u=1}^m\Big(dx+dy\,a_u+dz\,b_u\Big)I_u+$$
$$+\sum\limits_{r=m+1}^n\Big(dy\,a_r+dz\,b_r\Big)I_r
=\sum\limits_{u=1}^md\xi_u\,I_u+\sum\limits_{r=m+1}^ndT_r\,I_r\,,
$$
we have the following equality
$$
\zeta^{-1}d\zeta=\sum\limits_{u=1}^m\widetilde{A}_u\,d\xi_u\,I_u+
\sum\limits_{r=m+1}^n\widetilde{A}_{u_r}\,dT_r\,I_r+$$
\begin{equation}\label{obr-elem-0}
+\sum\limits_{s=m+1}^n\widetilde{A}_{s}\,d\xi_{u_s}\,I_s+
\sum\limits_{s=m+1}^n\sum\limits_{r=m+1}^n\widetilde{A}_{s}\,dT_r\,I_sI_r=:
\sum\limits_{k=1}^n\sigma_k\,I_k\,.
\end{equation}\vskip 2mm

Now, taking into account the denotation (\ref{obr-elem-0}) and the equality
(\ref{A__p}), we calculate:
$$\int\limits_{C_\zeta(0,R)}\sum\limits_{u=1}^m\sigma_u\,I_u=
\sum\limits_{u=1}^mI_u\int\limits_{C_u(\xi_u,R)}\frac{d\xi_u}{\xi_u}=
2\pi i\sum\limits_{u=1}^mI_u=2\pi i.
$$

Therefore,
\begin{equation}\label{obr-elem--1}
\lambda=2\pi i+\sum\limits_{k=m+1}^nI_k\int\limits_{C_\zeta(0,R)}\sigma_k\,.
\end{equation}
We note that form the relations (\ref{obr-elem--1}), (\ref{obr-elem}), and (\ref{A__p}) that
$\lambda$ is  an invertible element.

Thus, the equality (\ref{lambda-}) holds if and only if
\begin{equation}\label{obr-elem-1}
\int\limits_{C_\zeta(0,R)}\sigma_k=0\qquad \forall\;k=m+1,\ldots,n.
\end{equation}
  But,
for satisfying the equality (\ref{obr-elem-1}) the differential form $\sigma_k$ must be a total differential of some function.
We note that the property of being a total differential
is invariant under admissible transformations of coordinates
\cite[Theorem 2, p. 328]{Shabat}. In our situation, if we show that
$\sigma_k$ is a total differential of some function depend of the variables
$\frac{T_{m+1}}{\xi},\ldots,\frac{T_k}{\xi}$, then it means that
$\sigma_k$ is a total differential of some function depending on $x,y,z$.

\vskip2mm
\subsection{}
\vskip2mm

In this subsection we indicate a set of algebras in which the vectors
(\ref{e_1_e_2_e_3}) chosen arbitrarily and the equality (\ref{lambda-})
holds. We remind that an arbitrary commutative associative algebra, $\mathbb{A}_n^m$,
 with unit over the field of complex number $\mathbb{C}$ can be represented as
 $\mathbb{A}_n^m=S\oplus_s N$, where $S$ is $m$-dimensional semi-simple subalgebra
  and $N$ is $(n-m)$-dimensional nilpotent subalgebra (see Section~2).
\vskip2mm

 \theor\label{teo-pro-napivprostu-alg}
\textit{If $\mathbb{A}_n^m\equiv S$, then
the equality \em  (\ref{lambda-}) \em holds.}\vskip 1mm

The proof immediately follows from the conditions $\sigma_k\equiv 0$
 for $k=m+1,\ldots,n$ and (\ref{obr-elem--1}). This theorem is obtained in the paper \cite{Pl-Pukh-Analele}.
\vskip2mm

 \theor\label{teo-pro-nulyovy-nil'potentnu-pidalg}
\textit{If $\mathbb{A}_n^m=S\oplus_s N$ and $N$ is a zero nilpotent subalgebra, then
the equality \em  (\ref{lambda-}) \em holds.}\vskip 1mm

\textbf{Proof.} From the condition of theorem follows that in the
relations (\ref{A__p}) all $B_{k,p}=0$. Therefore, (\ref{A__p}) takes
the form
\begin{equation}\label{obr-elem-2}
\widetilde{A}_k=-\frac{T_k}{\xi_{u_k}^2},\qquad k=m+1,\ldots,n.
\end{equation}

Since $I_sI_r=0$ for $r,s=m+1,\ldots,n$, then form the denotation
(\ref{obr-elem-0}) and the identity (\ref{obr-elem-2}), we obtain
$$
\sigma_k=\frac{dT_k}{\xi_{u_k}}+\widetilde{A}_k\,d\xi_{u_k}=
\frac{dT_k}{\xi_{u_k}}-\frac{T_k}{\xi_{u_k}^2}\,d\xi_{u_k}=
d\left(\frac{T_k}{\xi_{u_k}}\right)=:d\tau_k\,,\quad k=m+1,\ldots,n.
$$
Under the transformation $(x,y,z)\rightarrow\tau_k$ the circle
$C_\zeta(0,R)$ maps into a closed smooth curve $\widetilde{C}$ (Jordan or not)
and the singularity $\xi_{u_k}=0$  maps on $\tau_k=\infty$.
Consequently, in an interior of the curve $\widetilde{C}$ does not exist
singular points. By the Cauchy theorem in complex plane \cite[p. 90]{Shabat},
we have:
$$\int\limits_{C_\zeta(0,R)}\sigma_k=\int\limits_{\widetilde{C}}
d\tau_k=0.
$$
So, the equality (\ref{lambda-}) is a consequence of the last relation and (\ref{obr-elem--1}).
\noindent The theorem is proved.

The Theorem \ref{teo-pro-nulyovy-nil'potentnu-pidalg}
implies the formula (\ref{form-Koshi--}) for monogenic functions
in the three-dimensional algebra $\mathbb{A}_2$ which investigated in the paper \cite{Pl-Pukh}.

Further we consider the case where $N$ in non-zero nilpotent subalgebra.
For this goal we establish an explicitly form of $\sigma_{m+1}, \sigma_{m+2}, \sigma_{m+3}$ and $\sigma_{m+4}$.

From the relation (\ref{obr-elem-0}) follows the equalities
\begin{equation}\label{sigma-k}
\begin{array}{c}
\displaystyle
\sigma_{m+1}=\frac{dT_{m+1}}{\xi_{u_{m+1}}}+\widetilde{A}_{m+1}\,d\xi_{u_{m+1}}\,,\vspace*{4mm}\\ \displaystyle
\sigma_{k}=\frac{dT_k}{\xi_{u_k}}+\widetilde{A}_k\,d\xi_{u_k}+
\sum\limits_{r,s=m+1}^{k-1}\widetilde{A}_r\,dT_s\Upsilon_{r,k}^s\,,\quad k=m+2,\ldots,n.\\
\end{array}
\end{equation}

Now, the equalities (\ref{A__p}) and (\ref{lem_3_2--}) implies the following
 equalities:
$$
\widetilde{A}_{m+1}=-\frac{T_{m+1}}{\xi_{u_{m+1}}^2}\,,\quad
\widetilde{A}_{m+2}=-\frac{T_{m+2}}{\xi_{u_{m+2}}^2}+\frac{T_{m+1}^2}
{\xi_{u_{m+2}}^3}\Upsilon_{m+1,m+2}^{m+1}\,,$$
\large

$$
\widetilde{A}_{m+3}=-\frac{T_{m+3}}{\xi_{u_{m+3}}^2}+\frac{T_{m+1}^2}
{\xi_{u_{m+3}}^3}\Upsilon_{m+1,m+3}^{m+1}+2\frac{T_{m+1}T_{m+2}}
{\xi_{u_{m+3}}^3}\Upsilon_{m+2,m+3}^{m+1}-$$
$$
-\frac{T_{m+1}^3}
{\xi_{u_{m+3}}^4}\Upsilon_{m+1,m+2}^{m+1}\Upsilon_{m+2,m+3}^{m+1}+
\frac{T_{m+2}^2}
{\xi_{u_{m+3}}^3}\Upsilon_{m+2,m+3}^{m+2}-\frac{T_{m+1}^2T_{m+2}}
{\xi_{u_{m+3}}^4}\Upsilon_{m+2,m+3}^{m+2}\Upsilon_{m+1,m+2}^{m+1}\,,
$$

$$
\widetilde{A}_{m+4}=-\frac{T_{m+4}}{\xi_{u_{m+4}}^2}+\frac{T_{m+1}^2}
{\xi_{u_{m+4}}^3}\Upsilon_{m+1,m+4}^{m+1}+2\frac{T_{m+1}T_{m+3}}
{\xi_{u_{m+4}}^3}\Upsilon_{m+3,m+4}^{m+1}+$$
$$
+2\frac{T_{m+1}T_{m+2}}
{\xi_{u_{m+4}}^3}\Upsilon_{m+2,m+4}^{m+1}+
2\frac{T_{m+2}T_{m+3}}
{\xi_{u_{m+4}}^3}\Upsilon_{m+3,m+4}^{m+2}+\frac{T_{m+2}^2}
{\xi_{u_{m+4}}^3}\Upsilon_{m+2,m+4}^{m+2}-$$
$$-\frac{T_{m+1}^3}{\xi_{u_{m+4}}^4}\Upsilon_{m+1,m+2}^{m+1}\Upsilon_{m+2,m+4}^{m+1}-
\frac{T_{m+1}^2T_{m+2}}
{\xi_{u_{m+4}}^4}\Upsilon_{m+1,m+2}^{m+1}\Upsilon_{m+2,m+4}^{m+2}-$$
$$
-\frac{T_{m+1}^2T_{m+3}}{\xi_{u_{m+4}}^4}\Upsilon_{m+1,m+2}^{m+1}\Upsilon_{m+3,m+4}^{m+2}+
\frac{T_{m+3}^2}{\xi_{u_{m+4}}^3}\Upsilon_{m+3,m+4}^{m+3}-
\frac{T_{m+1}^3}{\xi_{u_{m+4}}^4}\Upsilon_{m+1,m+3}^{m+1}\Upsilon_{m+3,m+4}^{m+1}-$$
$$
-\frac{T_{m+1}^2T_{m+2}}{\xi_{u_{m+4}}^4}\Upsilon_{m+1,m+3}^{m+1}\Upsilon_{m+3,m+4}^{m+2}-
\frac{T_{m+1}^2T_{m+3}}{\xi_{u_{m+4}}^4}\Upsilon_{m+1,m+3}^{m+1}\Upsilon_{m+3,m+4}^{m+3}-
$$
$$
-2\frac{T_{m+1}^2T_{m+2}}{\xi_{u_{m+4}}^4}\Upsilon_{m+2,m+3}^{m+1}\Upsilon_{m+3,m+4}^{m+1}-
2\frac{T_{m+1}T_{m+2}^2}{\xi_{u_{m+4}}^4}\Upsilon_{m+2,m+3}^{m+1}\Upsilon_{m+3,m+4}^{m+2}-$$
$$
-2\frac{T_{m+1}T_{m+2}T_{m+3}}{\xi_{u_{m+4}}^4}\Upsilon_{m+2,m+3}^{m+1}\Upsilon_{m+3,m+4}^{m+3}+
\frac{T_{m+1}^4}{\xi_{u_{m+4}}^5}\Upsilon_{m+1,m+2}^{m+1}\Upsilon_{m+2,m+3}^{m+1}\Upsilon_{m+3,m+4}^{m+1}+$$
$$
+\frac{T_{m+1}^3T_{m+2}}{\xi_{u_{m+4}}^5}\Upsilon_{m+1,m+2}^{m+1}\Upsilon_{m+2,m+3}^{m+1}\Upsilon_{m+3,m+4}^{m+2}+
\frac{T_{m+1}^3T_{m+3}}{\xi_{u_{m+4}}^5}\Upsilon_{m+1,m+2}^{m+1}\Upsilon_{m+2,m+3}^{m+1}\Upsilon_{m+3,m+4}^{m+3}-$$
$$
-\frac{T_{m+1}T_{m+2}^2}{\xi_{u_{m+4}}^4}\Upsilon_{m+2,m+3}^{m+2}\Upsilon_{m+3,m+4}^{m+1}-
\frac{T_{m+2}^3}{\xi_{u_{m+4}}^4}\Upsilon_{m+2,m+3}^{m+2}\Upsilon_{m+3,m+4}^{m+2}-$$
$$
-\frac{T_{m+2}^2T_{m+3}}{\xi_{u_{m+4}}^4}\Upsilon_{m+2,m+3}^{m+2}\Upsilon_{m+3,m+4}^{m+3}+
\frac{T_{m+1}^3T_{m+2}}{\xi_{u_{m+4}}^5}\Upsilon_{m+2,m+3}^{m+2}\Upsilon_{m+1,m+2}^{m+1}\Upsilon_{m+3,m+4}^{m+1}+$$
$$
+\frac{T_{m+1}^2T_{m+2}^2}{\xi_{u_{m+4}}^5}\Upsilon_{m+2,m+3}^{m+2}\Upsilon_{m+1,m+2}^{m+1}\Upsilon_{m+3,m+4}^{m+2}+
\frac{T_{m+1}^3T_{m+2}T_{m+3}}{\xi_{u_{m+4}}^5}\Upsilon_{m+2,m+3}^{m+2}\Upsilon_{m+1,m+2}^{m+1}\Upsilon_{m+3,m+4}^{m+3}\,.
$$
\vskip2mm

Finally, a consequence of the previous equalities and the relations
(\ref{sigma-k}) is the following differential representation of
 $\sigma_{m+1}, \sigma_{m+2}, \sigma_{m+3}$ and $\sigma_{m+4}$:
 \begin{equation}\label{sigma_m+1}
 \sigma_{m+1}=d\left(\frac{T_{m+1}}{\xi_{u_{m+1}}}\right),\quad
 \sigma_{m+2}=d\left(\frac{T_{m+2}}{\xi_{u_{m+2}}}-\frac{1}{2}\Upsilon_
 {m+1,m+2}^{m+1}
\frac{T_{m+1}^2}{\xi_{u_{m+2}}^2}\right),
 \end{equation}
 $$\sigma_{m+3}=d\Biggr(\frac{T_{m+3}}{\xi_{u_{m+3}}}-\frac{1}{2}\Upsilon_
 {m+1,m+3}^{m+1}
\frac{T_{m+1}^2}{\xi_{u_{m+3}}^2}-\Upsilon_{m+2,m+3}^{m+1}
\frac{T_{m+1}T_{m+2}}{\xi_{u_{m+3}}^2}-$$
$$-\frac{1}{2}\Upsilon_{m+2,m+3}^{m+2}
\frac{T_{m+2}^2}{\xi_{u_{m+3}}^2}+\frac{1}{3}\Upsilon_{m+1,m+2}^{m+1}\Upsilon_{m+2,m+3}^{m+1}
\frac{T_{m+1}^3}{\xi_{u_{m+3}}^3}\Biggr)+$$
\begin{equation}\label{sigma_m+3}
+ \Upsilon_{m+1,m+2}^{m+1}\Upsilon_{m+2,m+3}^{m+2}\,\sigma_{m+3}^{(1)}\,,
\end{equation}\large
$$
\sigma_{m+4}=d\Biggr(\frac{T_{m+4}}{\xi_{u_{m+4}}}-\frac{1}{2}\Upsilon_
 {m+1,m+4}^{m+1}
\frac{T_{m+1}^2}{\xi_{u_{m+4}}^2}-\Upsilon_{m+3,m+4}^{m+1}
\frac{T_{m+1}T_{m+3}}{\xi_{u_{m+4}}^2}-\vspace*{2mm} $$
$$
-\Upsilon_{m+2,m+4}^{m+1}
\frac{T_{m+1}T_{m+2}}{\xi_{u_{m+4}}^2}
-\frac{1}{2}\Upsilon_{m+2,m+4}^{m+2}
\frac{T_{m+2}^2}{\xi_{u_{m+4}}^2}-\Upsilon_{m+3,m+4}^{m+2}
\frac{T_{m+2}T_{m+3}}{\xi_{u_{m+3}}^2}+\vspace*{2mm} $$
$$
+\frac{1}{3}\Upsilon_{m+1,m+2}^{m+1}\Upsilon_{m+2,m+4}^{m+1}
\frac{T_{m+1}^3}{\xi_{u_{m+4}}^3}-
\frac{1}{2}\Upsilon_{m+3,m+4}^{m+3}
\frac{T_{m+3}^2}{\xi_{u_{m+4}}^2}+
\frac{1}{3}\Upsilon_{m+1,m+3}^{m+1}\Upsilon_{m+3,m+4}^{m+1}
\frac{T_{m+1}^3}{\xi_{u_{m+4}}^3}-\vspace*{2mm} $$
$$
-\frac{1}{4}\Upsilon_{m+1,m+2}^{m+1}\Upsilon_{m+2,m+3}^{m+1}
\Upsilon_{m+3,m+4}^{m+1}
\frac{T_{m+1}^4}{\xi_{u_{m+4}}^4}+
\frac{1}{3}\Upsilon_{m+2,m+3}^{m+2}\Upsilon_{m+3,m+4}^{m+2}
\frac{T_{m+2}^3}{\xi_{u_{m+4}}^3}\Biggr)+
\vspace*{2mm}$$
$$
+\Upsilon_{m+1,m+2}^{m+1}\Upsilon_{m+2,m+4}^{m+2}\,\sigma_{m+4}^{(1,2)}+
\Upsilon_{m+1,m+2}^{m+1}\Upsilon_{m+3,m+4}^{m+2}\,\sigma_{m+4}^{(2,2)}+
\Upsilon_{m+1,m+3}^{m+1}\Upsilon_{m+3,m+4}^{m+2}\,\sigma_{m+4}^{(3,2)}+\vspace*{2mm} $$
$$
+\Upsilon_{m+3,m+4}^{m+3}\Upsilon_{m+1,m+3}^{m+1}\,\sigma_{m+4}^{(4,3)}+
\Upsilon_{m+2,m+3}^{m+1}\Upsilon_{m+3,m+4}^{m+1}\,\sigma_{m+4}^{(5,1)}+
\Upsilon_{m+2,m+3}^{m+1}\Upsilon_{m+3,m+4}^{m+2}\,\sigma_{m+4}^{(6,2)}+
\vspace*{2mm} $$
$$
+\Upsilon_{m+2,m+3}^{m+1}\Upsilon_{m+3,m+4}^{m+3}\,\sigma_{m+4}^{(7,3)}-
\Upsilon_{m+1,m+2}^{m+1}\Upsilon_{m+2,m+3}^{m+1}
\Upsilon_{m+3,m+4}^{m+2}\,\sigma_{m+4}^{(8,2)}-$$
$$
-\Upsilon_{m+1,m+2}^{m+1}\Upsilon_{m+2,m+3}^{m+1}
\Upsilon_{m+3,m+4}^{m+3}\,\sigma_{m+4}^{(9,3)}+
\Upsilon_{m+2,m+3}^{m+2}\Upsilon_{m+3,m+4}^{m+1}\,\sigma_{m+4}^{(10,1)}+\vspace*{2mm}
$$
$$+\Upsilon_{m+2,m+3}^{m+2}\Upsilon_{m+3,m+4}^{m+3}\,\sigma_{m+4}^{(11,3)}
-\Upsilon_{m+2,m+3}^{m+2}\Upsilon_{m+1,m+2}^{m+1}
\Upsilon_{m+3,m+4}^{m+1}\,\sigma_{m+4}^{(12,1)}-\vspace*{2mm}$$
\begin{equation}\label{sigma-m+4}
\begin{array}{c}-\Upsilon_{m+2,m+3}^{m+2}\Upsilon_{m+1,m+2}^{m+1}
\Upsilon_{m+3,m+4}^{m+2}\,\sigma_{m+4}^{(13,2)}
-\Upsilon_{m+2,m+3}^{m+2}\Upsilon_{m+1,m+2}^{m+1}
\Upsilon_{m+3,m+4}^{m+3}\,\sigma_{m+4}^{(14,3)}\,,\\
\end{array}
\end{equation}
where
\begin{equation}\label{sigma-m+3-1}
\sigma_{m+3}^{(1)}:=\frac{T_{m+1}^2}{\xi^3_{u_{m+3}}}\left(dT_{m+2}-
 \frac{T_{m+2}}{\xi_{u_{m+3}}}\,d\xi_{u_{m+3}}\right),
 \end{equation}
and $\sigma_{m+4}^{(\ell,r)}$, \,$\ell=1,2,\ldots,14$ are determined by the following
relations:
\begin{equation}\label{sigma-m+4-1}
\sigma_{m+4}^{(\ell,r)}:=\left\{
\begin{array}{lcl}
\frac{T_{m+1}^2}{\xi_{u_{m+4}}^3}\,g(r)& \quad\mbox{for}\quad&\ell=1,2,3,4,\vspace*{2mm}\\
\frac{2\,T_{m+1}T_{m+2}}{\xi_{u_{m+4}}^3}\,g(r)& \mbox{for}&\ell=5,6,7,\vspace*{2mm}\\
\frac{T_{m+1}^3}{\xi_{u_{m+4}}^4}\,g(r)& \mbox{for}&\ell=8,9,\vspace*{2mm}\\
\frac{T_{m+2}^2}{\xi_{u_{m+4}}^3}\,g(r)& \mbox{for}&\ell=10,11,\vspace*{2mm}\\
\frac{T_{m+1}^2T_{m+2}}{\xi_{u_{m+4}}^4}\,g(r)& \mbox{for}&\ell=12,13,14,\vspace*{2mm}\\
\end{array}
\right.
\end{equation}
where $g(r):=dT_{m+r}-\frac{T_{m+r}}{\xi_{u_{m+4}}}\,d\xi_{u_{m+4}}$.
\vskip3mm

\theor\label{teo-pro-nil'potentnu-pidalg-menshe-rivne-3}
\textit{If $\mathbb{A}_n^m=S\oplus_s N$ and $\dim_{\mathbb{C}}N\leq3$,
 then the equality \em (\ref{lambda-}) \em holds.}\vskip 1mm

\textbf{Proof.} From the equality (\ref{sigma_m+1}) for $\sigma_{m+1}$, we have
$$
\sigma_{m+1}=d\left(\frac{T_{m+1}}{\xi_{u_{m+1}}}\right)=:d\tau_{m+1}\,.$$
Now, the identity
$\int_{C_\zeta(0,R)}\sigma_{m+1}=0
$
proved as in Theorem \ref{teo-pro-nulyovy-nil'potentnu-pidalg}.

Consider $\sigma_{m+2}$ from the equality (\ref{sigma_m+1}), which is a total
differential of the certain function depending on the variables
$\frac{T_{m+1}}{\xi_{u_{m+2}}}\,,
\frac{T_{m+2}}{\xi_{u_{m+2}}}$.
Under the transformation $(x,y,z)\rightarrow\left(\frac{T_{m+1}}{\xi_{u_{m+2}}}\,,
\frac{T_{m+2}}{\xi_{u_{m+2}}}\right)$ 
 the circle
$C_\zeta(0,R)$ maps into a closed smooth curve $\widetilde{\widetilde{C}}$ (Jordan or not)
and the singularity $\xi_{u_{m+2}}=0$  maps on $\infty$.
Consequently, in an interior of the curve $\widetilde{\widetilde{C}}$ does not exist
singular points. Then by the Cauchy theorem in the space $\mathbb{C}^2$
 \cite[p. 334]{Shabat},
we have:
$$\int\limits_{C_\zeta(0,R)}\sigma_{m+2}(x,y,z)=\int\limits_{\widetilde{\widetilde{C}}}
\sigma_{m+2}\left(\frac{T_{m+1}}{\xi_{u_{m+2}}}\,,
\frac{T_{m+2}}{\xi_{u_{m+2}}}\right)=0.
$$

Finally, we prove the equality (\ref{obr-elem-1}) for $k=m+3$. In the paper
\cite{Burde_de_Graaf} is described all commutative associative nilpotent algebras
over the field $\mathbb{C}$ of dimensional $1,2,3$.
From results of the paper \cite{Burde_de_Graaf} (Table~1) immediately follows
 that for all mentioned algebras the relation
$\Upsilon_{m+1,m+2}^{m+1}\Upsilon_{m+2,m+3}^{m+2}=0$ is always satisfied.
Therefore, the equality (\ref{sigma_m+3}) implies that under the conditions of
theorem $\sigma_{m+3}$ is always a total
differential of the certain function depending on the variables
$\frac{T_{m+1}}{\xi_{u_{m+3}}}\,,
\frac{T_{m+2}}{\xi_{u_{m+3}}}\,,\frac{T_{m+3}}{\xi_{u_{m+3}}}$.

Now as before, under the transformation
$$(x,y,z)\rightarrow\left(\frac{T_{m+1}}{\xi_{u_{m+3}}}\,,
\frac{T_{m+2}}{\xi_{u_{m+3}}}\,,\frac{T_{m+3}}{\xi_{u_{m+3}}}\right)$$ 
 the circle
$C_\zeta(0,R)$ maps into a closed smooth curve $\widehat{C}$ (Jordan or not)
and the singularity $\xi_{u_{m+3}}=0$  maps on $\infty$.
Hence, in an interior of the curve $\widehat{C}$ does not exist
singular points. Then by the Cauchy theorem in the space $\mathbb{C}^3$
 \cite[p. 334]{Shabat},
we have:
$$\int\limits_{C_\zeta(0,R)}\sigma_{m+3}(x,y,z)=\int\limits_{\widehat{C}}
\sigma_{m+3}\left(\frac{T_{m+1}}{\xi_{u_{m+3}}}\,,
\frac{T_{m+2}}{\xi_{u_{m+3}}}\,,\frac{T_{m+3}}{\xi_{u_{m+3}}}\right)=0.
$$
So, the equality (\ref{lambda-}) is a consequence of the last relation and (\ref{obr-elem--1}).

\noindent The theorem is proved.

Let us note that from the Theorem \ref{teo-pro-nil'potentnu-pidalg-menshe-rivne-3}
follows the formula (\ref{form-Koshi--}) for monogenic functions
in the three-dimensional algebra $\mathbb{A}_3$ (see \cite{Pl-Shp3}) and
 in the three-dimensional algebra $\mathbb{A}_2$ which considered in the paper
 \cite{Pl-Pukh}.
 \vskip2mm

\theor\label{teo-pro-nil'potentnu-pidalg-rivne-4}
\textit{Let $\mathbb{A}_n^m=S\oplus_s N$ and $\dim_{\mathbb{C}}N=4$.
Then the equality \em (\ref{lambda-}) \em holds if the following relations
satisfied}\large
\begin{equation}\label{v-teor-pro-sigma-m+4}
\begin{array}{l}
\Upsilon_{m+1,m+2}^{m+1}\Upsilon_{m+2,m+3}^{m+2}=
\Upsilon_{m+1,m+2}^{m+1}\Upsilon_{m+2,m+4}^{m+2}=
\Upsilon_{m+1,m+3}^{m+1}\Upsilon_{m+3,m+4}^{m+2}=\vspace*{3mm} \\ 
=\Upsilon_{m+3,m+4}^{m+3}\Upsilon_{m+1,m+3}^{m+1}=
\Upsilon_{m+2,m+3}^{m+1}\Upsilon_{m+3,m+4}^{m+1}=
\Upsilon_{m+2,m+3}^{m+1}\Upsilon_{m+3,m+4}^{m+2}=\vspace*{3mm} \\
=\Upsilon_{m+2,m+3}^{m+1}\Upsilon_{m+3,m+4}^{m+3}=
\Upsilon_{m+1,m+2}^{m+1}\Upsilon_{m+2,m+3}^{m+1}
\Upsilon_{m+3,m+4}^{m+2}=\vspace*{3mm} \\
=\Upsilon_{m+1,m+2}^{m+1}\Upsilon_{m+2,m+3}^{m+1}
\Upsilon_{m+3,m+4}^{m+3}
=\Upsilon_{m+2,m+3}^{m+2}\Upsilon_{m+3,m+4}^{m+1}=
\Upsilon_{m+2,m+3}^{m+2}\Upsilon_{m+3,m+4}^{m+3}=\vspace*{3mm} \\
=\Upsilon_{m+2,m+3}^{m+2}\Upsilon_{m+1,m+2}^{m+1}
\Upsilon_{m+3,m+4}^{m+1}=
\Upsilon_{m+2,m+3}^{m+2}\Upsilon_{m+1,m+2}^{m+1}
\Upsilon_{m+3,m+4}^{m+2}=\vspace*{3mm} \\
=\Upsilon_{m+2,m+3}^{m+2}\Upsilon_{m+1,m+2}^{m+1}
\Upsilon_{m+3,m+4}^{m+3}=0.\\
\end{array}
\end{equation}
\vskip 1mm

\textbf{Proof.} From the equalities (\ref{sigma_m+3}) and (\ref{sigma-m+4}) it obvious
that under conditions (\ref{v-teor-pro-sigma-m+4}) expressions for $\sigma_{m+3}$ and
$\sigma_{m+4}$ are total differentials. Further proof is similar to proof of the
 Theorem
\ref{teo-pro-nil'potentnu-pidalg-menshe-rivne-3}. The theorem is proved.

Further we consider some examples of algebras, which satisfy the relations (\ref{v-teor-pro-sigma-m+4}).\vskip2mm

\textbf{Examples.}
\begin{itemize}
  \item Consider the algebra with the basis $\{I_1:=1,I_2,I_3,I_4,I_5\}$ and
  multiplication rules:
  $$I_2^2=I_3\,,\,\,I_2\,I_4=I_5$$ and other products are zeros (for nilpotent
  subalgebra see \cite{Martin}, Table 21, algebra $\mathcal{J}_{69}$ and
  \cite{Burde_Fialowski}, page 590, algebra $A_{1,4}$).
  \item Consider the algebra with the basis $\{I_1:=1,I_2,I_3,I_4,I_5\}$ and
   multiplication rules:
  $$I_2^2=I_3$$ and other products are zeros (for nilpotent
  subalgebra see \cite{Burde_Fialowski}, page 590, algebra $A_{1,2}\oplus A_{0,1}^2$).
  \item The algebra with the basis $\{I_1:=1,I_2,I_3,I_4,I_5\}$ and
   multiplication rules:
  $$I_2^2=I_3\,,\,\,I_4^2=I_5$$ and other products are zeros (for nilpotent
  subalgebra see \cite{Burde_Fialowski}, page 590, algebra $A_{1,2}\oplus A_{1,2}$).
  \item The algebra with the basis $\{I_1:=1,I_2,I_3,I_4,I_5\}$ and
  multiplication rules:
  $$I_2^2=I_3\,,\,\,I_2\,I_3=I_4$$ and other products are zeros (for nilpotent
  subalgebra see \cite{Martin}, Table 21, algebra $\mathcal{J}_{71}$).
\end{itemize}

Now we consider an example of algebra, which does not satisfy the relations (\ref{v-teor-pro-sigma-m+4}). Moreover, we choose the vectors $e_1,e_2,e_3$ of the
form (\ref{e_1_e_2_e_3}) such that the equality
 (\ref{lambda-}) is not true. \vskip2mm

\textbf{Example.} \vskip2mm

Consider the algebra $\mathbb{A}_5$ with the basis
$\{1,\rho,\rho^2,\rho^3,\rho^4\}$, where $\rho^5=0$
(see \cite{Pl-Shp-Algeria} and \cite{Rosculet-55}, par. 11). Here $n=5$, $m=1$.
 It is obvious that $\Upsilon_{2,3}^2\Upsilon_{3,4}^3=1$ and
the relations (\ref{v-teor-pro-sigma-m+4}) are not true.
Consider the vectors:
$$e_1=1, \quad e_2=i+\rho^2+\rho^4,\quad e_3=(1-i)\rho+
\left(\frac{1}{4}-\frac{3}{4}\,i\right)\rho^3,
$$
which are linearly independent over $\mathbb{R}$ and satisfy the equality
$$e_1^2+e_2^2+e_3^2=0.$$
Let $\zeta=xe_1+ye_2+ze_3$. In the algebra $\mathbb{A}_5$ for given $\zeta$,
we have \, $$\xi_{u_2}=\xi_{u_3}=\xi_{u_4}=\xi_{u_5}=x+iy=:\xi\,.$$

The inverse element $\zeta^{-1}$ is of the form (\ref{obr-elem}), where
$$\widetilde{A}_0=\frac{1}{\xi}\,,\quad \widetilde{A}_1=\frac{z(i-1)}{\xi^2}\,,
\quad \widetilde{A}_2=-\frac{y}{\xi^2}+\frac{z^2(1-i)^2}{\xi^3}\,,
$$

$$\widetilde{A}_3=\frac{1}{4}\frac{z(3i-1)}{\xi^2}+\frac{2yz(1-i)}{\xi^3}-
\frac{z^3(1-i)^3}{\xi^4}\,,
$$

$$\widetilde{A}_4=-\frac{y}{\xi^2}+\frac{y^2+\frac{1}{2}z^2(1-i)(1-3i)}{\xi^3}-
\frac{3yz^2(1-i)^2}{\xi^4}+\frac{z^4(1-i)^4}{\xi^5}\,.
$$

Let us set
\begin{equation}\label{kolo-c}
C_\zeta(0,R):=\{\zeta=xe_1+ye_2\in E_3: x^2+y^2=R^2\}.
\end{equation}
On the circle of integration (\ref{kolo-c}), we obtain:
\begin{equation}\label{A-_-}
\widetilde{A}_0=\frac{1}{\xi}\,,\quad \widetilde{A}_1=\widetilde{A}_3=0,\quad
\widetilde{A}_2=-\frac{y}{\xi^2}\,,\quad\widetilde{A}_4=-\frac{y}{\xi^2}+
\frac{y^2}{\xi^3}\,.
\end{equation}
As a consequence of the equations (\ref{sigma-k}), (\ref{A-_-}) on the
circle (\ref{kolo-c}) we obtain the following expression
$$\sigma_5=\left(\frac{1}{\xi}-\frac{y}{\xi^2}\right)dy+
\left(-\frac{y}{\xi^2}+\frac{y^2}{\xi^3}\right)d\xi.
$$

It is easy to calculate that
$$\int\limits_{C_\zeta(0,R)}\sigma_5=\frac{\pi i}{2}
$$
and
$$\int\limits_{C_\zeta(0,R)}\sigma_1=\int\limits_{|\xi|=R}\frac{d\xi}{\xi}=2\pi i,
\qquad
\int\limits_{C_\zeta(0,R)}\sigma_k=0,\qquad k=2,3,4.
$$
Hence, in this example
$$\lambda=\int\limits_{C_\zeta(0,R)}\zeta^{-1}d\zeta=2\pi i+\frac{\pi i}{2}\rho^4.
$$

\vskip2mm
\subsection{}
\vskip2mm

In this subsection we indicate sufficient conditions on a choose of
the vectors (\ref{e_1_e_2_e_3}) for which the equality (\ref{lambda-}) is true.
Let the algebra $\mathbb{A}_n^m$ be represented as $\mathbb{A}_n^m=S\oplus_s N$. Let us
note that the condition  $\zeta\in E_3\subset S$ means that in
the decomposition (\ref{e_1_e_2_e_3}) \,
 $a_k=b_k=0$ for all $k=m+1,\ldots,n$.
\vskip2mm

 \theor\label{teo-pro-napivprostu-alg-}
\textit{If $\mathbb{A}_n^m=S\oplus_s N$ and $\zeta\in E_3\subset S$,
 then the equality \em (\ref{lambda-}) \em holds.}\vskip 1mm

\textbf{Proof.} Since $\zeta\in S$, then $T_k=0$ for $k=m+1,\ldots,n$ (see denotation
(\ref{lem_3_2-})). From (\ref{lem_3_2--}), (\ref{A__p}) follows that $\widetilde{A}_k=0$,
 and now from (\ref{sigma-k}) follows that $\sigma_k=0$ for $k=m+1,\ldots,n$. The equality (\ref{lambda-}) is a consequence of the equality $\sigma_k=0$ and the relation (\ref{obr-elem--1}).
\noindent The theorem is proved.

Let us note that by essentially the Theorem \ref{teo-pro-napivprostu-alg-}
generalizes the Theorem 3 of the paper \cite{Shpakiv-Kuzm-Analele}.

Now we consider a case where $\zeta\notin S$. If $\mathbb{A}_n^m=S\oplus_s N$ and $\dim_{\mathbb{C}}N\leq3$,  then by Theorem  \ref{teo-pro-nil'potentnu-pidalg-menshe-rivne-3}
 the equality (\ref{lambda-}) holds for any $\zeta\in E_3$.
 \vskip2mm

\theor\label{teo-pro-nil'potentnu-pidalg-rivne-4-dzeta}
\textit{Let $\mathbb{A}_n^m=S\oplus_s N$ and $\dim_{\mathbb{C}}N=4$.
Then the equality \em (\ref{lambda-}) \em holds if the following two conditions
satisfied:
  \begin{enumerate}
    \item $a_{m+1}=b_{m+1}=0$;
    \item at least one of the relations $a_{m+2}=b_{m+2}=0$ or $a_{m+3}=b_{m+3}=0$ are true.
  \end{enumerate}}

\textbf{Proof.} It follows from the condition of theorem that $T_{m+1}=0$ and
at least one of the equalities $T_{m+2}=0$ or $T_{m+3}=0$ are true.
To prove (\ref{lambda-}) it is need to prove the equality (\ref{obr-elem--1}) for $k=m+1,\ldots,m+4$. The equality (\ref{obr-elem--1}) is proved in Theorem \ref{teo-pro-nil'potentnu-pidalg-menshe-rivne-3} for $k=m+1,m+2$.
Under the condition $T_{m+1}=0$ from (\ref{sigma-m+3-1}), we have $\sigma_{m+3}^{(1)}=0$.
Since now $\sigma_{m+3}$ is a total differential, then similar to proof of Theorem
\ref{teo-pro-nil'potentnu-pidalg-menshe-rivne-3} can be proved the
equality (\ref{obr-elem--1}) for $k=m+3$.

Moreover, under the conditions of
theorem from the denotation (\ref{sigma-m+4-1}) follows the equalities $\sigma_{m+4}^{(\ell,r)}=0$
for all $\ell=1,\ldots,14$. Therefore, $\sigma_{m+4}$ is a total differential, then similar to proof of Theorem
\ref{teo-pro-nil'potentnu-pidalg-menshe-rivne-3} can be proved the
equality (\ref{obr-elem--1}) for $k=m+4$.

\noindent The theorem is proved.

\vskip10mm

\vskip7mm
Vitalii Shpakivskyi

 Department of Complex Analysis and Potential Theory

 Institute of Mathematics of the National Academy of Sciences of
Ukraine,

 3, Tereshchenkivs'ka st.

01601 Kyiv-4

 UKRAINE

 http://www.imath.kiev.ua/\~{}complex/

 \ e-mail: shpakivskyi@mail.ru,\, shpakivskyi@imath.kiev.ua

\end{document}